\documentstyle[11pt]{article}
\begin{document}
\setlength{\baselineskip}{15pt}
\title{Quasipolynomial Generalization of Lotka-Volterra Mappings}
\author{ Benito Hern\'{a}ndez--Bermejo$^{\:(1)}$  \and L\'{e}on Brenig$^{\:(2)}$}
\date{}

\maketitle
\begin{center}
{\em \mbox{} \\
Service de Physique Th\'{e}orique et Math\'{e}matique. Universit\'{e} Libre de Bruxelles. \\
Campus Plaine -- CP 231. Boulevard du Triomphe, B-1050 Brussels, Belgium.}
\end{center}

\mbox{}

\begin{abstract}
In the last years it has been shown that Lotka-Volterra mappings constitute a valuable tool from both the theoretical and the applied points of view, with developments in very diverse fields such as Physics, Population Dynamics, Chemistry and Economy. The purpose of this work is to demonstrate that many of the most important ideas and algebraic methods that constitute the basis of the quasipolynomial formalism (originally conceived for the analysis of ordinary differential equations) can be extended into the mapping domain. The extension of the formalism into the discrete-time context is remarkable as far as the quasipolynomial methodology had never been shown to be applicable beyond the differential case. It will be demonstrated that Lotka-Volterra mappings play a central role in the quasipolynomial formalism for the discrete-time case. Moreover, the extension of the formalism into the 
discrete-time domain allows a significant generalization of Lotka-Volterra mappings as well as a whole transfer of algebraic methods into the discrete-time context. The result is a novel and more general conceptual framework for the understanding of Lotka-Volterra mappings, as well as a new range of possibilities that becomes open not only for the theoretical analysis of 
Lotka-Volterra mappings and their generalizations, but also for the development of new applications. 
\end{abstract}

\mbox{}



\mbox{}

Short title: Generalization of LV mappings.

\mbox{}

\mbox{}

\mbox{}

\mbox{}

\mbox{}

\mbox{}

\mbox{}

\mbox{}

\mbox{}

\mbox{}

\noindent $^{(1)}$ E-mail: bhernand@ulb.ac.be

\noindent $^{(2)}$ Corresponding author. Fax: (+ 00 32 2) 650 58 24. E-mail: lbrenig@ulb.ac.be

\pagebreak

\newtheorem{df}{Definition}
\newtheorem{rk}{Remark}
\newtheorem{th}{Theorem}
\newtheorem{lm}{Lemma}
\newtheorem{co}{Corollary}
\newtheorem{cj}{Conjecture}

\begin{flushleft}
{\bf 1. Introduction}
\end{flushleft}

In recent years, the Quasipolynomial (QP in what follows) formalism has provided a tool for the obtainment of useful information about very general systems of ODEs in an algebraic way. Historically, it was introduced independently by several authors \cite{bro}-\cite{pym}. Since then, the formalism has experienced a remarkable degree of extension and generalization leading to the establishment of a plethora of results, both mathematical (involving integrability properties \cite{bre}-\cite{bre2}, reduction and simplification techniques 
\cite{pla206,pla253}, normal form analysis \cite{slb1}, approximate construction of solutions \cite{slb2}, determination of first integrals \cite{bvmp1}-\cite{atb2}, stability and construction of Lyapunov functions \cite{bvmaa1}-\cite{gat1}, development of numerical algorithms \cite{ryv1}, analysis of algebraic properties \cite{bvl1}) and applied (Hamiltonian systems \cite{bvmp1,bvmaa1,pla234}, neural networks \cite{mlvb1}, system modelling 
\cite{bv3,bv2}, population dynamics \cite{pym,bvmp1,bvmaa1}, connectionism \cite{rbv1}, cosmology \cite{frgbf1}, etc).

Specifically, Lotka-Volterra (LV in what follows) differential systems constitute a particular but prominent case of the QP equations, therefore playing a canonical role in the formalism. This has allowed the generalization and transfer of many properties previously known for 
Lotka-Volterra systems to more general situations \cite{slb1}-\cite{gat1}, 
\cite{bvl1}, \cite{bv3}-\cite{rbv1}. 

Simultaneously to the development of the ideas and methods of the QP formalism, and in a completely independent way, during the last years the theory of LV mappings has also been developed. As it was the case for their differential counterpart, LV mappings were first introduced in a biological context \cite{may1,mo1}. In this domain they have proven to be richer from the theoretical point of view and more realistic from the modelling perspective than other classical possibilities such as the logistic equation \cite{hhj1,lst1}. Subsequently, LV mappings have become of interest for the analysis of dynamical properties and chaotic behavior in different domains, mainly Population Dynamics \cite{kr1}-\cite{lw1}, Physics \cite{raj1}-\cite{ur2}, Chemistry \cite{gb1} and Economy \cite{bd1,doh1}.

The purpose of this work is to demonstrate that many of the most important ideas and algebraic methods that constitute the basis of the QP formalism can be translated into the discrete, i.e. the mapping domain. In this sense, the establishment of the mathematical basis of this generalization is the central goal of the paper. The extension of the formalism into the discrete context is remarkable as far as the QP methodology had never been shown to be applicable beyond the differential case. It will be demonstrated that LV mappings play a canonical role in the QP formalism for the mapping case, analogous to the one of  LV systems in the differential QP scenario. However, the extension of the formalism into the discrete domain allows also an important generalization of LV mappings and, in addition, an entire transfer of algebraic methods into the novel discrete-time context. The outcome is a new and more general conceptual framework for the understanding of LV mappings as well as a whole range of possibilities that becomes open both for the theoretical analysis of LV mappings and their generalizations, and for new applications in system modelling. 

The structure of the article is as follows. In Section 2 the foundations of the QP formalism for mappings is presented. Section 3 is devoted to the analysis of the relationship between 
the differential and the discrete QP equations, from the discretization point of view. The work concludes in Section 4 with some final remarks in which the future implications of the results are outlined.

\mbox{}

\mbox{}

\begin{flushleft}
{\bf 2. QP mappings}
\end{flushleft}

\begin{flushleft}
{\em 2.1 Fundamentals}
\end{flushleft}

We start by defining the main subject of this work:

\mbox{}

\begin{df}
\label{df1}
{\rm QP mappings are those of the form
\begin{equation}
\label{qpm}
	x_i(p+1)=x_i(p) \exp \left( \lambda _i + \sum _{j=1}^m A_{ij} \prod _{k=1}^n 
	[x_k(p)]^{B_{jk}} \right)
	\:\: , \:\:\:\:\: i =1, \ldots , n
\end{equation}
where {\em (i)} $m$ is an integer not necessarily equal to $n$; {\em (ii)} index $p$ is an 
integer denoting the discrete time; {\em (iii)} variables $x_i(p)$ are assumed to be positive for $i=1, \ldots ,n$ and for every $p$; and {\em (iv)} $A=(A_{ij})$, 
$B=(B_{ij})$ and $\lambda = ( \lambda _i)$ are real matrices of dimensions $n \times m$, $m \times n$ and $n \times 1$, respectively. 
}
\end{df}

\mbox{}

Therefore we see that the well-known LV mappings 
\begin{equation}
\label{lvm}
	x_i(p+1)=x_i(p) \exp \left( \lambda _i + \sum _{j=1}^n A_{ij} x_j(p) \right)
	\:\: , \:\:\:\:\: i =1, \ldots , n
\end{equation}
are a particular case of QP mapping, namely the one corresponding to $m=n$ and $B$ the $n \times n$ identity matrix.

For future developments it is convenient to introduce an additional matrix, denoted by $M$, which is of dimension $n \times (m+1)$ and is defined as:
\begin{equation}
	M = \left( \begin{array}{cccc}
		\lambda _1 & A_{11} & \ldots  & A_{1m} \\
		\lambda _2 & A_{21} & \ldots  & A_{2m} \\
		\vdots     & \vdots & \mbox{} & \vdots \\
		\lambda _n & A_{n1} & \ldots  & A_{nm} 
	\end{array} \right)
\end{equation}

Notice also that every QP mapping is completely specified by its associated matrices $A$, $B$ and $\lambda$ or, equivalently, $M$ and $B$. The algebraic manipulation of the matrices instead of that of the whole system is one of the characteristic features of the QP formalism. Before starting the description of the formalism, it is convenient to give a simple result which complements Definition \ref{df1} and will be evoked later:

\mbox{}

\begin{lm}
\label{lm1}
{\rm For every QP mapping, if there is one integer $p$ such that $x_i(p)>0$ for all $i=1, \ldots ,n$, then $x_j(q)>0$ for every integer $q$ and for all $j=1, \ldots ,n$.
}
\end{lm}

\pagebreak
\begin{flushleft}
{\em Proof.}
\end{flushleft}

The proof is straightforward.$\Box$

\mbox{}

\begin{co}
\label{vco1}
{\rm The positive orthant is an invariant set for every QP mapping.
}
\end{co}

\mbox{}

Although this issue will be considered in more detail in Section 4, notice that this corollary is relevant as far as in many applications ---such as Population Dynamics--- positivity is a necessary property. From the point of view of the QP formalism, it will always be considered that the dynamics takes place in the interior of the positive orthant, according to Definition \ref{df1}, Lemma \ref{lm1} and Corollary \ref{vco1}. 

We can now start the description of the QP formalism for mappings (\ref{qpm}). 

\mbox{}

\begin{df}
\label{df2}
{\rm The terms 
\begin{equation}
\label{qm}
	\prod _{k=1}^n [x_k(p)]^{B_{jk}} \: , \:\:\: j=1, \ldots ,m
\end{equation}
appearing in the argument of the exponential in equation (\ref{qpm}) shall be named quasimonomials (QMs from now on).
}
\end{df}

\mbox{}

\begin{df}
\label{df3}
{\rm A quasimonomial (QM) transformation of a set of positive variables $\{ x_1(p) , \ldots ,$ 
$x_n(p) \}$, which are discrete in time (i.e. $p$ is an integer) is defined as
\begin{equation}
\label{qmt}
	x_i(p)= \prod_{j=1}^n [y_j(p)]^{C_{ij}} \:\: , \:\:\: i = 1, \ldots ,n
\end{equation}
where $C$ is a real and invertible $n \times n$ matrix, and $\{ y_1(p) , \ldots , y_n(p) \}$ is another set of positive discrete-time variables.
}
\end{df}

\mbox{}

The form-invariance of QP mappings after a QM transformation is one of the cornerstones of the formalism, as the following two results show:

\mbox{}

\begin{th}
\label{th1}
{\rm Consider a $n$-dimensional QP mapping of matrices $A$, $B$, $\lambda$ and $M$. After a QM transformation of matrix $C$, the result is another $n$-dimensional QP mapping of matrices $A'$, $B'$, $\lambda '$ and $M'$, where:
\begin{equation}
\label{mtqmt}
	A' = C^{-1} \cdot A \: , \:\:\: 
	B' = B \cdot C \: , \:\:\: 
	\lambda ' = C^{-1} \cdot \lambda \: , \:\:\: 
	M' = C^{-1} \cdot M 
\end{equation}
}
\end{th}

\pagebreak
\begin{flushleft}
{\em Proof.}
\end{flushleft}

When we apply transformation (\ref{qmt}) on mapping (\ref{qpm}) we find:
\begin{equation}
	\prod_{r=1}^n [y_r(p+1)]^{C_{ir}}=\prod_{s=1}^n [y_s(p)]^{C_{ir}} \exp \left( \lambda _i 	+ \sum_{j=1}^m A_{ij} \prod _{s=1}^n [y_s(p)]^{B'_{js}} \right) \:\: , \:\:\: 
	B' = B \cdot C
\end{equation}
Since all the variables are positive, we can take logarithms on both sides and obtain after some algebra:
\begin{equation}
	y_i(p+1)=y_i(p) \exp \left( \lambda '_i + \sum _{j=1}^m A'_{ij} \prod _{k=1}^n 
	[y_k(p)]^{B'_{jk}} \right) 	\:\: , \:\:\:\:\: i =1, \ldots , n
\end{equation}
with $A' = C^{-1} \cdot A$ , $B' = B \cdot C$ and $\lambda ' = C^{-1} \cdot \lambda$, as announced.$\Box$

\mbox{}

However, form-invariance alone is not sufficient in order to have a complete and consistent set of rules for the transformation of QP mappings by means of QM transformations. A second theorem complements the previous one:

\mbox{}

\begin{th}
\label{th2}
{\rm Every QM transformation relating two QP mappings is a topological conjugacy.
}
\end{th}

\begin{flushleft}
{\em Proof.}
\end{flushleft}

QM transformations are continuous and bijective in the positive orthant, as needed in order to have a topological conjugacy. In addition, consider two $n$-dimensional QP mappings of matrices $\{A,B, \lambda \}$ and $\{A',B', \lambda ' \}$ which are related by a QM transformation of matrix $C$, i.e. $A' = C^{-1} \cdot A$ , $B' = B \cdot C$ and $\lambda ' = 
C^{-1} \cdot \lambda$. We have:
\begin{eqnarray*}
	F(x) & = & \left\{ x_i \exp \left( \lambda _i + \sum _{j=1}^m A_{ij} \prod _{k=1}^n 
	x_k^{B_{jk}} \right) , 	\:\: i=1, \ldots ,n \right\} \\
	G(y) & = & \left\{ y_i \exp \left( \lambda '_i + \sum _{j=1}^m A'_{ij} \prod _{k=1}^n 
	y_k^{B'_{jk}} \right) , 	\:\: i=1, \ldots ,n \right\} \\
	\varphi(x) & = & \left\{ \prod_{j=1}^n x_j^{C_{ij}^{-1}}, \:\: i=1, \ldots ,n \right\}
\end{eqnarray*}
After some algebra it is possible to verify that:
\begin{equation}
	\varphi_i(F(x))=G_i( \varphi (x))= \left( \prod_{j=1}^n x_j^{C^{-1}_{ij}} \right) 
	\exp \left( \lambda '_i + \sum _{l=1}^m A'_{il} \prod _{k=1}^n x_k^{B_{lk}} \right) 
\: , \:\: i=1, \ldots ,n
\end{equation}
This set of equalities demonstrates that the QM transformation is a topological 
conjugacy.$\Box$

\mbox{}

Consequently, we not only have a form invariance between QP systems related by a QM transformation, but actually a complete dynamical equivalence. This dynamical equivalence motivates the following definition:

\mbox{}

\begin{df}
\label{df4}
{\rm Consider a QP mapping. The set of all QP mappings related to it by means of QM transformations shall be termed its equivalence class. 
}
\end{df}

\mbox{}

This definition is consistent, since the relation defined among QP systems connected by QM transformations is actually an equivalence one. Notice that:

\mbox{}

\begin{co}
\label{dco1}
{\rm The matrix product $B \cdot M$ is invariant for every equivalence class.
}
\end{co}

\mbox{}

At this point we have a first (still incomplete) framework that allows the establishment of some results which are interesting as an illustration of how the formalism works, and also because they are to be used as the basis for further developments. 

\mbox{}

\begin{flushleft}
{\em 2.2 Reduction to the nonredundant form}
\end{flushleft}

Quasimonomial transformations can be used to reduce an arbitrary QP mapping to what shall be termed the non-redundant form of the mapping. Let us start by defining it:

\mbox{}

\begin{df}
\label{df5}
{\rm A QP mapping (\ref{qpm}) is said to be non-redundant if $m \geq n$ and the ranks of both 
$B$ and $M$ are maximal (i.e. Rank($B$) = Rank($M$) = $n$). Otherwise the QP mapping shall be termed redundant.
}
\end{df}

\mbox{}

That this definition is natural shall become clear along this subsection. The best way to see it is by demonstrating that every redundant mapping can be always reduced to a non-redundant one. The reduction proceeds in three steps:

\mbox{}

\begin{flushleft}
{\em 2.2.1 Step 1: Reduction to the case $m \geq n$}
\end{flushleft}

Assume a QP mapping for which $m < n$. Let rank($B$) $=r$, with $r \leq m$. 
Since $B$ is an $m \times n$ matrix, dim\{Ker($B$)\}$= (n-r) > 0$. 
Consequently, we can perform a QM transformation of matrix
\begin{equation}
   \label{mmn}
   C = \left( \begin{array}{c} 
            I_r \\ \mbox{} \\ O_{(n-r) \times r}
       \end{array} \right| 
       \left. \begin{array}{ccc}
            \mbox{}        & \mbox{} & \mbox{}     \\
            \phi ^{(r+1)}  & \ldots  & \phi ^{(n)} \\
            \mbox{}        & \mbox{} & \mbox{} 
       \end{array} \right)
\end{equation}
where $I_r$ is the $r \times r$ identity matrix, and $\{ \phi ^{(r+1)}, 
\ldots , \phi ^{(n)} \}$ is a set of column vectors which constitute a basis 
of Ker($B$). When applied, the quasimonomial transformation given by matrix 
(\ref{mmn}) leads to:
\begin{equation}
   B' = B \cdot C = 
      \left( \begin{array}{ccc}
        B_{11} & \ldots  & B_{1r} \\
        \vdots & \mbox{} & \vdots \\
        B_{m1} & \ldots  & B_{mr}
      \end{array} \right| 
      \left. \begin{array}{ccc}
          0    & \ldots  & 0      \\
        \vdots & \mbox{} & \vdots \\
          0    & \ldots  & 0
      \end{array} \right)
\end{equation}
Consequently, in the transformed mapping the ($n-r$) final variables are not 
present in the arguments of the exponentials: 
\begin{equation}
	y_i(p+1)=y_i(p) \exp \left( \lambda '_i + \sum _{j=1}^m A'_{ij} \prod _{k=1}^r 
	[y_k(p)]^{B_{jk}} \right) \:\: , \:\:\:\:\: i =1, \ldots , n
\end{equation}
This implies that, in fact, we are led to an $r$-dimensional QP mapping for $\{ y_1 , \ldots , 
y_r \}$ plus ($n-r$) relationships for $\{ y_{r+1}, \ldots ,y_n \}$, which are completely determined by the values of $\{ y_1 , \ldots , y_r \}$. Therefore $\{ y_{r+1} , \ldots , y_n 
\}$ become non-essential variables that can be decoupled from the mapping, which is thus reduced to dimension $n'=r$. Consequently we have $m \geq n'$ as stated. 

In the $r$-dimensional reduced mapping we have the following $B$ matrix:
\begin{equation}
   \hat{B}' = 
      \left( \begin{array}{ccc}
        B_{11} & \ldots  & B_{1r} \\
        \vdots & \mbox{} & \vdots \\
        B_{m1} & \ldots  & B_{mr}
      \end{array} \right)
\end{equation}
Notice that $\hat{B}'$ has been obtained from the matrix $B$ of the initial QP mapping. Actually, $\hat{B}'$ is the result of deleting the $(n-r)$ top right columns of $B$. Therefore it might still happen that different QMs, i.e. rows in $B$, could become degenerated, namely equal, in $\hat{B}'$, as a result of the column deletion. However, notice that $\hat{B}'$ is the matrix denoting the system QMs, and therefore cannot have equal rows, since two equal rows denote one and the same QM. Nevertheless, since $\hat{B}'$ is a $m \times r$ matrix, with $m \geq r$, and Rank($\hat{B}'$)$=r$, even if such degeneration of some of the rows takes place, the rank value guarantees that the number of independent rows will be al least $r$. Therefore, if we have to cancel out some equal rows of $\hat{B}'$, the result would be a QP mapping with $n'=r$ and a value of $m'$ such that $m \geq m' \geq r=n'$, and consequently $m' \geq n'$ as expected. 

\mbox{}

\begin{flushleft}
{\em 2.2.2 Step 2: Reduction to a mapping with maximum Rank($B$)}
\end{flushleft}

Assume a QP mapping for which $m \geq n$ and Rank($B$) is not maximun. Let Rank($B$)$=r<n$. Then the mapping can be transformed into another one for which the rank of matrix $B$ is maximum. The method is quite similar to the one used in Subsection 2.2.1.: Since the rank of 
$B$ is not maximum, then ($n-r$) variables can be decoupled by means of a QM transformation of matrix
\begin{equation}
   \label{mmn2}
   C = \left( \begin{array}{c} 
            I_r \\ \mbox{} \\ O_{(n-r) \times r}
       \end{array} \right| 
       \left. \begin{array}{ccc}
            \mbox{}        & \mbox{} & \mbox{}     \\
            \phi ^{(r+1)}  & \ldots  & \phi ^{(n)} \\
            \mbox{}        & \mbox{} & \mbox{} 
       \end{array} \right)
\end{equation}
where the notation is the same as the one employed in (\ref{mmn}). We thus arrive at a QP mapping of matrix:
\begin{equation}
   B' = B \cdot C = 
      \left( \begin{array}{ccc}
        B_{11} & \ldots  & B_{1r} \\
        \vdots & \mbox{} & \vdots \\
        B_{m1} & \ldots  & B_{mr}
      \end{array} \right| 
      \left. \begin{array}{ccc}
          0    & \ldots  & 0      \\
        \vdots & \mbox{} & \vdots \\
          0    & \ldots  & 0
      \end{array} \right)
\end{equation}
In this case, this is possible because dim\{Ker($B$)\}$=(n-r)$. Then we are led to a reduced QP mapping of $n'=r$ variables and $B$ matrix given by
\begin{equation}
   \hat{B}' = 
      \left( \begin{array}{ccc}
        B_{11} & \ldots  & B_{1r} \\
        \vdots & \mbox{} & \vdots \\
        B_{m1} & \ldots  & B_{mr}
      \end{array} \right)
\end{equation}
where Rank($\hat{B}'$)$=r$ is maximum. Again, in the case that some degeneracies of the QMs might appear as a consequence of the column deletion, an argument similar to the one employed in Subsection 2.2.1 shows that we always arrive to a QP mapping of $n'$ variables and $m'$ QMs, with $m \geq m' \geq r =n'$ as we wanted. Moreover, let $\hat{B}'_*$ be the $B$ matrix of such mapping, obtained after the suppression in $\hat{B}'$ of the degenerated QMs. Then the rank of $\hat{B}'_*$ is clearly maximum, i.e. equal to $n'=r$. Therefore the reduction to a QP mapping with $m' \geq n'$ and maximum rank of its $B$ matrix (given by $\hat{B}'_*$) is achieved. We can thus proceed into the last step of the reduction to the non-redundant form.

\mbox{}

\begin{flushleft}
{\em 2.2.3 Step 3: Reduction to a mapping with maximum Rank($M$)}
\end{flushleft}

Assume a QP mapping with $m \geq n$. It is convenient to begin by stating the following result:

\mbox{}

\begin{th}
\label{th3}
{\rm Consider a $n$-dimensional QP mapping for which $m \geq n$. If Rank($M$)$=r<n$, then there exists a QM transformation leading to a QP mapping of matrix
\begin{equation}
\label{mcc}
	M'= \left( \begin{array}{c} 
		M'_{r \times (m+1)} \\ O_{(n-r) \times (m+1)} \end{array} \right)
\end{equation}
where $O$ is the null matrix and the subindexes of the submatrices denote their sizes.
}
\end{th}

\begin{flushleft}
{\em Proof.}
\end{flushleft}

If the reduction is possible there exists an invertible $n$-dimensional 
matrix $C$ such that $M' = D \cdot M$, where $D = C^{-1}$ and $M'$ 
complies to form (\ref{mcc}). We shall show that $C$ exists by constructing 
it. We have the matrix identity:
\begin{equation}
       \left( \begin{array}{c}
               M'_{r \times (m+1)} \\
               O_{(n-r) \times (m+1)}
          \end{array} \right) \;\: = \;\: 
       \left( \begin{array}{c}
               D_{r \times n} \\
               D_{(n-r) \times n}
          \end{array} \right) \; \cdot \; M 
\end{equation}
Here we write $D$ as composed of two submatrices. This equation implies that 
the $(m+1)$ column vectors of matrix $M$ must belong to Ker\{$D_{(n-r) \times n}$\}, 
which is of dimension $r$. Since Rank($M$) $=r$, there 
exists a basis set of column vectors of the form $\{ {\bf v}_1 , \ldots , 
{\bf v}_r \}$ that spans all the column vectors of $M$. We can complete this 
basis to a basis of $I \!\! R^n$, namely $\{ {\bf v}_1 , \ldots , {\bf v}_r, 
{\bf w}_{r+1}, \ldots , {\bf w}_n \}$. In this basis set, the projection 
matrix  
\begin{equation}
       \Pi ' = 
       \left( \begin{array}{cc}
               O_{r \times r} & \mbox{} \\
               \mbox{} & I_{(n-r) \times (n-r)}
          \end{array} \right)
\end{equation}
has as kernel the desired subspace generated by $\{ {\bf v}_1 , \ldots , 
{\bf v}_r \}$. We only have to change to the canonical basis to obtain a 
matrix $\Pi = P^{-1} \cdot \Pi ' \cdot P$, where $P$ is a $n \times n$ 
invertible matrix. But since Rank($\Pi$) $=$ ($n-r$), the system of equations 
of the kernel of $\Pi$
\begin{equation}
      \Pi \cdot {\bf x} = {\bf 0}
\end{equation}
possesses $r$ redundant equations. If we supress them, the result is a 
system of the form
\begin{equation}
      \Pi _R \cdot {\bf x} = {\bf 0} \; .
\end{equation}
But $\Pi _R$ complies to our requirements for $D_{(n-r) \times n}$. We can 
therefore take $ D_{(n-r) \times n} = \Pi _R$. The whole matrix $D$ is thus 
formed by completing the submatrix $D_{r \times n}$ with arbitrary entries, 
in such a way that the resulting matrix $D$ is invertible. Consequently, the 
searched quasimonomial transformation exists and the proof is complete.$\Box$

\mbox{}

Therefore, after the QM transformation of Theorem \ref{th3} we arrive at the QP mapping:
\begin{eqnarray}
\label{qpth3}
	y_i(p+1) & = & y_i(p) \exp \left( \lambda '_i + \sum _{j=1}^m A'_{ij} \prod _{k=1}^n 
	[y_k(p)]^{B'_{jk}} \right) \:\: , \:\:\:\:\: i =1, \ldots , r \\
 	y_j(p+1) & = & y_j(p) \:\: , \:\:\:\:\: j = r+1, \ldots , n 
\end{eqnarray}
We can then decouple the last $(n-r)$ variables, which are in fact constants, and write
\begin{equation}
	y_i(p+1) = y_i(p) \exp \left( \hat{ \lambda }'_i + \sum _{j=1}^m \hat{A}'_{ij} 
	\prod _{k=	1}^r [y_k(p)]^{ \hat{B}'_{jk}} \right) \:\: , \:\:\:\:\: i =1, \ldots , r 
\end{equation}
where $\hat{\lambda}'_i=\lambda '_i$, $i=1, \ldots ,r$; $\hat{B}'_{jk}=B'_{jk}$, 
$j=1, \ldots m$, $k=1, \ldots ,r$; and: 
\begin{equation}
   \hat{A}' = 
      \left( \begin{array}{ccc}
        A'_{11} q_1 & \ldots  & A'_{1m}q_m \\
        \vdots & \mbox{} & \vdots \\
        A'_{r1} q_1 & \ldots  & A'_{rm}q_m \\
      \end{array} \right) \:\:\: , \:\:\:\:\:\: q_j=\prod_{k=r+1}^n [y_k(0)]^{B'_{jk}}
\end{equation}
Notice that the $M$ matrix of this QP mapping, which is given by $\hat{M}'=( \hat{\lambda} ' \mid \hat{A}' )$, is now of maximal rank, i.e. Rank($\hat{M}'$)$=r$. Consequently, we have arrived at an $r$-dimensional QP mapping with maximal rank of $\hat{M}'$. 

In addition, if the $B$ matrix of the original QP mapping is of maximal rank, then $B'=B \cdot C$ is also of maximal rank. As it was the case in steps 1 and 2, the decoupling of ($n-r$) variables can lead to a degeneracy of QMs in matrix $\hat{B}'$. Since every column of matrix 
$\hat{A}'$ corresponds to one QM, if such degeneracy takes place then both $\hat{B}'$ and 
$\hat{A}'$ have to be rewritten in order to regroup the coefficients of the same QM. For instance, if we have $m=n=3$, $r=2$ and matrices:
\begin{equation}
	B' = \left( \begin{array}{ccc} 1 & 1 & 1 \\ 1 & 1 & 0 \\ 1 & 0 & 0 \end{array} \right) 
	\:\: , \:\:\:\:\:
	M' = \left( \begin{array}{cccc} 
	\lambda '_1 & A'_{11} & A'_{12} & A'_{13} \\ 
	\lambda '_2 & A'_{21} & A'_{22} & A'_{23} \\ 
	0 & 0 & 0 & 0 
	\end{array} \right)
\end{equation}
We have decoupled the third variable, so we must suppress the third column of $B'$, and the result is:
\begin{equation}
	B' \rightarrow
	\hat{B}' = 
	\left( \begin{array}{cc} 1 & 1 \\ 1 & 1 \\ 1 & 0 \end{array} \right) 
	\:\: , \:\:\:\:\:
	M' \rightarrow
	\hat{M}' = \left( \begin{array}{cccc} 
	\lambda '_1 & A'_{11}q_1 & A'_{12}q_2 & A'_{13}q_3 \\ 
	\lambda '_2 & A'_{21}q_1 & A'_{22}q_2 & A'_{23}q_3 
	\end{array} \right)
\end{equation}
But then, the first and second rows of $\hat{B}'$ are the same, thus representing one and the same QM. This is by definition a redundant way of representing a QP mapping. Therefore, the optimal way to express the mapping matrices is:
\begin{equation}
	\hat{B}'_* = 
	\left( \begin{array}{cc} 1 & 1 \\ 1 & 0 \end{array} \right) 
	\:\: , \:\:\:\:\:
	\hat{M}'_* = \left( \begin{array}{ccc} 
	\lambda '_1 & (A'_{11}q_1 + A'_{12}q_2) & A'_{13}q_3 \\ 
	\lambda '_2 & (A'_{21}q_1 + A'_{22}q_2) & A'_{23}q_3 
	\end{array} \right)
\end{equation}
After the example it should be simpler to understand the general problem. Normally, after the application of the reduction described in Theorem \ref{th3}, we arrive at a QP mapping of matrices $\hat{B}'$ and $\hat{M}'$, as indicated above. If $B$ is of maximal rank, then both $\hat{B}'$ and $\hat{M}'$ are of maximal rank, namely $r$. If there are no coincidences in the rows of $\hat{B}'$, the reduction is completed. If this is not the case, we have to perform a further simplification of matrices $\hat{B}'$ and $\hat{M}'$, so we arrive to the final matrices $\hat{B}'_*$ and $\hat{M}'_*$. Then it is not difficult to see that:
\begin{itemize}
\item $\hat{B}'_*$ is a $m' \times n'$ matrix, with $n'=r$ and $m \geq m' \geq r=n'$. In addition, Rank($\hat{B}'_*$)$=r$ is maximum.
\item $\hat{M}'_*$ is a $n' \times m'$ matrix, with $n'=r$ and $m \geq m' \geq r=n'$. However,    Rank($\hat{M}'_*$) needs not be maximum. Therefore we have two cases:
\begin{enumerate}
	\item If Rank($\hat{M}'_*$) is maximum (i.e. equal to $r$) the reduction is complete.
	\item If Rank($\hat{M}'_*$) is not maximum, then Theorem \ref{th3} can be applied again. In 	this way we can apply Step 3 in an iterative way until the reduction is completed.
\end{enumerate}
\end{itemize}

Therefore, after the application of the three previous steps, we arrive at the non-redundant form of the QP mapping, in which $m \geq n$ and the ranks of $B$ and $M$ are maximum, i.e. equal to the number of variables $n$. 

The meaning of Theorem \ref{th3} can also be understood in an alternative way:

\mbox{}

\begin{th}
\label{th4}
{\rm Consider a QP mapping for which $m \geq n$. If Rank($M$)$=r<n$, then there are $(n-r)$  functionally independent constants of motion of quasimonomial form.
}
\end{th}

\begin{flushleft}
{\em Proof.}
\end{flushleft}

If we apply Theorem \ref{th3} we arrive at mapping (\ref{qpth3}), in which $y_j= constant$ for all $j=r+1, \ldots ,n$. From the definition of QM transformation we have:
\begin{equation}
y_j(p)= \prod _{k=1}^n x_k^{C_{jk}^{-1}}(p) = y_j(0) \:\: , \:\:\:\: j=r+1, \ldots ,n
\end{equation}
The functional independence is just a consequence of the linear independence of the rows of 
$C^{-1}$, which is assumed by hypothesis.$\Box$

\mbox{}

Therefore the reduction made in Step 3 essentially consists of a decoupling of these constants of motion by means of a suitable QM transformation.

\mbox{}
\begin{rk}
\label{rk2}
{\rm In the following, it is always assumed that we are working with the QP mappings in their non-redundant form, unless otherwise stated.
}
\end{rk}

\mbox{}

We are now in position to complete the description of the QP formalism for discrete-time systems.

\mbox{}

\begin{flushleft}
{\em 2.3 Further transformations and equivalence to the LV canonical form}
\end{flushleft}

In this subsection a central result of the formalism will be demonstrated, i.e. the equivalence of every QP mapping to an associated or canonical LV mapping, to which it can be related by means of a topological conjugacy. To prove it, some new transformations called embeddings shall be introduced. In this way, the demonstration of the result shall lead to the completion of our description of the foundations of the QP formalism for mappings. 

In that respect, two cases must be distinguished:

\mbox{}

\begin{flushleft}
{\em 2.3.1 Case 1: $m=n$}
\end{flushleft}

Consider an arbitrary QP mapping for which $m=n$. As usual, let $\{ M, B \}$ be the matrices of the mapping, and let Rank($B$)$=n$. Consider a QM transformation of matrix $C=B^{-1}$, which is correctly defined because $B$ is by hypothesis invertible. From (\ref{mtqmt}), the result is another QP mapping of matrices
\begin{equation}
	M' = B \cdot M \:\: , \:\:\:\: B' = I_{n \times n}
\end{equation}
where $I$ denotes the identity matrix. Since $B$ is the identity matrix, what we have is a LV mapping, as stated. Recall that, according to Theorem \ref{th2}, the QM transformation is a topological conjugacy, and consequently the original and the LV mappings are conjugate to each other. Another interesting remark at this moment is that, according to Corollary \ref{dco1}, $M'$ is the class invariant. In this sense, the LV mapping can be considered as the canonical representative of the QP class of equivalence. 

\mbox{}

\begin{flushleft}
{\em 2.3.2 Case 2: $m>n$}
\end{flushleft}

Consider now the complementary case $m>n$. Let a QP mapping of matrices $\{ M,B \}$ such that Rank($B$)$=n$. Consider now the following mapping, obtained from the initial one just by adding $(m-n)$ trivial variables of constant value 1:
\begin{eqnarray}
\label{embqp}
	x_i(p+1) & = & x_i(p) \exp \left( \lambda _i + \sum _{j=1}^m A_{ij} \prod _{k=1}^n 
	[x_k(p)]^{B_{jk}} \right)
	\:\: , \:\:\:\:\: i =1, \ldots , n \\
	x_j(p+1) & = & x_j(p) \:\: , \:\:\:\:\:  x_j(0)=1 \:\: ; \:\: j = n+1, \ldots , m
	\label{embvar}
\end{eqnarray}
This is also a QP mapping to which the following matrices can be assigned:
\begin{equation}
\label{emb}
	\tilde{M} = \left( \begin{array}{c} M \\ O_{(m-n) \times (m+1)} \end{array} \right) 
	\:\: , \:\:\:\: 
	\tilde{B} = ( B \mid \bar{B}_{m \times (m-n)})
\end{equation}
In (\ref{emb}), the subindexes of the submatrices denote their sizes, $O$ is the null submatrix and $\bar{B}$ is an arbitrary submatrix chosen in such a way that Rank($\tilde{B}$)$=m$, which is always possible. 

\mbox{}

\begin{df}
\label{df6}
{\rm Let a $n$-dimensional QP mapping for which $m > n$ and Rank($B$)$=n$. The construction of an associated $m$-dimensional QP mapping according to the rule (\ref{emb}) shall be termed an embedding.
}
\end{df}

\mbox{}

We have the following important fact:

\mbox{}

\begin{lm}
\label{lm2}
{\rm Consider a QP mapping of variables $\{ x_1, \ldots , x_n \}$ for which 
$m>n$ and Rank($B$)$=n$. In addition, consider its embedded mapping of variables $\{ x_1, \ldots , x_m \}$. Then the initial mapping is conjugate to the embedded mapping in the level set $\{ x_j=1 \: , \:\: j=n+1, \ldots , m \}$.
}
\end{lm}

\begin{flushleft}
{\em Proof.}
\end{flushleft}

The proof is straightforward.$\Box$

\mbox{}

Therefore, an embedding amounts to an increase in the dimension of the mapping, but the dynamics is not affected. However, the embedded mapping has been reduced to the $m=n$ case. Moreover, the embedding is defined in such a way that matrix $\tilde{B}$ is invertible. Consequently, the problem has been reduced to Case 1. If we now perform a QM transformation of matrix $C = 
\tilde{B}^{-1}$ over the embedded system, the result is a $m$-dimensional LV mapping of matrices:
\begin{equation}
	\tilde{M}' = \tilde{B} \cdot \tilde{M} = B \cdot M \:\: , \:\:\:\: 
	\tilde{B}' = I_{m \times m}
\end{equation}
And again, the LV mapping coefficients are given by the class invariant $B \cdot M$ of the initial mapping. Notice, in addition, that now Rank($\tilde{M}'$) is not maximum due to the zero submatrix introduced in $\tilde{M}$ in the embedding process, as indicated in (\ref{emb}). In particular, it is simple to see that Rank($\tilde{M}'$)$=$Rank($M$). Consequently, Theorem \ref{th4} is applicable to the LV mapping, and thus a set of quasimonomial and functionally independent constants of motion is present in it. Actually, it should be noted that the embedding operation is essentially the inverse of the reduction procedure developed in Subsection 2.2.3.

These considerations allow to summarize the previous results in the following:

\mbox{}

\begin{th}
\label{th5}
{\rm Consider a QP mapping of matrices $B$ and $M$, with $m \geq n$ and such that 
Rank($B$)$=n$. Then:
\begin{description}
\item[{\em i)\/}] If $m=n$ then it is conjugate to a LV mapping of matrix $B \cdot M$.
\item[{\em ii)\/}] If $m>n$ then it is conjugate to a LV mapping of matrix $B \cdot M$ in the level set of $(m-n)$ functionally independent constants of motion of quasimonomial form present in the LV mapping.
\end{description}
}
\end{th}

\begin{flushleft}
{\em Proof.}
\end{flushleft}

Proof for {\em (i)} is complete in Subsection 2.3.1. For part {\em (ii)\/} all the ingredients were presented in Subsection 2.3.2, and to complete the proof it is only necessary to realize that the set of constants of motion is actually the transformed of the level set of variables $\{ x_j(p) = constant =  x_j(0)=1 \:\ , \:\: j = n+1, \ldots , m \}$ given in (\ref{embvar}) after the QM transformation of matrix $C = \tilde{B}^{-1}$. This completes the demonstration.$\Box$

\mbox{}

Although in both cases the LV mapping is given by the class invariant matrix $B \cdot M$ of the initial QP mapping, in the case $m=n$ the LV mapping can be considered as the canonical representative of the equivalence class (since the LV mapping actually belongs to the class). On the contrary, in the case $m>n$ the LV mapping does not belong to the class, since it has dimension $m$. However, following the terminology of the QP formalism for differential equations, in both cases we can consider the LV mapping as a unique canonical form associated to every QP mapping. The fact that the LV mapping is unique and that there is a topological dynamical equivalence between it and the initial QP mapping is one of the cornerstones of the formalism, both from a theoretical point of view as well as for many applications.

We now have the tools for the establishment of an important result completing the foundations of the formalism.

\mbox{}

\begin{flushleft}
{\em 2.4 Characterization of the equivalence classes}
\end{flushleft}

The previous developments regarding the reduction to the non-redundant form and the LV embedding can be used to establish the result that completes the description of the equivalence classes given in Subsection 2.1:

\mbox{}

\begin{th}
\label{th6}
{\rm Consider two $n$-dimensional QP mappings characterized by matrices $\{ M, B \}$ and 
$\{ M', B' \}$, respectively, and both having the same number $m$ of QMs. Assume, in addition, that $m \geq n$ and that Rank($B$) $=$ Rank($B'$) $=n$. Then the two  mappings belong to the same class of equivalence if and only if  $B \cdot M = B' \cdot M'$.
}
\end{th}

\begin{flushleft}
{\em Proof.}
\end{flushleft}

We only have to demonstrate that there exists a unique, invertible $n \times n$ matrix $C$ such that the QM transformation of matrix $C$ connects both mappings.

For the case $m = n$ it can be checked straightforwardly that there is only 
one possibility, $C = B^{-1} \cdot B'$, which is obviously invertible. This 
choice of $C$ leads to $A' = C^{-1} \cdot A$ and $\lambda ' = C^{-1} \cdot 
\lambda$, and then $C$ complies to all the requirements of the class. 

The case $m > n$ can be reduced to the previous one with $m$ variables and $m$ quasimonomials by means of an embedding. In the embedded class of equivalence, there exists a $m \times m$
matrix $\tilde C$ that relates both sets of variables. But since the $(m-n)$
new variables added in the embedding are constants of value 1, this shows 
that the structure of $\tilde C$ is as follows:
\begin{equation}
  \tilde C = \left( \begin{array}{cc}
              C_{ n \times n} & C_{n \times (m-n) } \\
              C_{ (m-n) \times n} & C_{ (m-n)\times (m-n) } 
             \end{array} \right) \;\: ,
\end{equation}
where $C_{n \times n}$ relates the original sets of variables. This 
demonstrates the existence of $C$. The uniqueness of $C$ holds from the fact 
that both $B$ and $B'$ are of maximum rank. This implies that the linear 
system $B \cdot C = B'$ possesses a single solution in $C$: Suppose this is 
not the case. Then we have $C \neq C'$ such that $B \cdot C = B'$ and $B 
\cdot C' = B'$. In particular, $B \cdot (C -C') = B \cdot D = 
O_{m \times n}$, and the only solution to this system is $D = O_{n \times n}$ 
or $C = C'$, in contradiction with the original assumption. Finally, $C$ must 
be invertible since $B \cdot C = B'$.$\Box$

\mbox{}

Therefore, the result given in Corollary \ref{dco1} is now completed. However, it has been necessary to first develop the formalism before stating Theorem \ref{th6}. 

In order to completely set the foundations of the QP formalism for discrete-time systems, we shall now consider its relationship with the traditional, differential QP formalism.

\mbox{}

\mbox{}

\begin{flushleft}
{\bf 3. QP mappings as a commutative discretization of QP differential equations}
\end{flushleft}

In this section we shall demonstrate that QP mappings can be seen as a discretization of differential QP systems. Of course, many possible discretizations can be defined for a differential system. However, we shall demonstrate that QP mappings can be regarded as a very special discretization, namely the one preserving the form invariance properties of the differential equations. This has important consequences that shall be discussed in the next section. 

Before discussing this subject, we shall recall here some basic facts about the differential QP formalism for ODEs that will be necessary in what is to follow in order to make the article self-contained. However, this reminder shall be reduced to a very minimum, and we refer the reader to the bibliography cited in the Introduction for further details. 

\mbox{}

\begin{flushleft}
{\em 3.1 Some basic facts about QP systems of ODEs}
\end{flushleft}

The starting point of the differential QP formalism is the definition of QP equations:
\begin{equation}
   \frac{d{x}_i}{dt} = \dot{x}_i = x_{i} \left( \lambda _{i} + \sum_{j=1}^{m}A_{ij}\prod_{k=
      1}^{n}x_{k}^{B_{jk}} \right) \: , \;\:\;\: i = 1 \ldots n 
   \label{eq:glv}
\end{equation}
where $n$ and $m$ are positive integers, and $A$, $B$ and $\lambda$ are 
$n \times m$, $m \times n$ and $n \times 1$ real matrices, respectively. The definition of the combined matrix of coefficients $M = ( \lambda \mid A)$ is also used, where $M$ is  $n \times (m+1)$. It is always assumed that the $x_i$ are real and positive. As in the mapping case, $n$ denotes the number of variables, and $m$ the number of QMs
\begin{equation}
   \prod_{k=1}^{n}x_{k}^{B_{jk}} , \;\:\;\: j = 1 \ldots m 
\end{equation}
System (\ref{eq:glv}) is formally invariant under QM transformations
\begin{equation}
   x_{i} = \prod_{j=1}^{n} y _{j}^{C_{ij}} , \;\:\;\: i=1,\ldots ,n
   \label{bec}
\end{equation}
for any invertible real matrix $C$. After (\ref{bec}), the result is another QP system for which matrices $B, A, \lambda$ and $M$ change to 
\begin{equation}
B' = B \cdot C \;\:, \:\;\: A' = C^{-1} \cdot A \;\: , \:\;\: 
\lambda ' = C^{-1} \cdot \lambda \;\: , \:\;\: M' = C^{-1} \cdot M 
\end{equation}
These transformation rules define classes of equivalence, composed by all those QP systems related through QM transformations (\ref{bec}). The matrix product $B \cdot M$ is thus an invariant of every class of equivalence. Under the previous assumptions, the solutions of all the members of a class are topologically equivalent. 

It can also be demonstrated that every $n$-dimensional QP differential system of matrices 
$\{ A,B, \lambda \}$ with $m \geq n$ and Rank($B$) $=n$ can be transformed into a 
$m$-dimensional LV canonical form:
\begin{equation}
   \dot{z}_i = z_{i} \left( \tilde{\lambda }'_{i} + \sum_{j=1}^{m} \tilde{A}'_{ij}z_j \right) 
	\: , \;\:\;\: i = 1 \ldots m 
   \label{eq:lv}
\end{equation}
with $\tilde{\lambda }'=B \cdot \lambda$ and $\tilde{A}'=B \cdot A$. Actually, the dynamics of the original QP system is topologically equivalent to that of the LV system in a 
$n$-dimensional submanifold of the phase space of the latter, given by a level set of $(m-n)$ QM first integrals.

\mbox{}

\begin{flushleft}
{\em 3.2 QP discretization of QP differential systems}
\end{flushleft}

We shall now consider the discretization procedure leading from QP differential systems to QP mappings. In order to better appreciate the properties of the discretization, it shall be compared with the standard Euler discretization. Our starting point is the QP differential system that is now written in the form:
\begin{equation}
   \label{glvdis}
	\dot{x}_i = x_{i} \left( \lambda ^*_{i} + \sum_{j=1}^{m}A^*_{ij}
	\prod_{k=1}^{n}x_{k}^{B_{jk}} \right) , \;\:\;\: i = 1 \ldots n 
\end{equation}
In (\ref{glvdis}) the coefficients $\lambda ^*_{i}$ and $A^*_{ij}$ have an asterisk superscript in order to denote that their physical dimension is the inverse of time. In other words, they can be considered as coefficients of interaction per unit time. Assume now that we wish to discretize the system in time steps of length $\epsilon \ll 1$. One possibility is, for instance, the standard Euler discretization of the derivative:
\begin{equation}
	\frac{x_i(t+ \varepsilon ) - x_i(t)}{\varepsilon} \approx \dot{x}_i(t) = 
	x_{i}(t) \left( \lambda ^*_{i} + \sum_{j=1}^{m}A^*_{ij}\prod_{k=1}^{n}x_{k}^{B_{jk}}(t) 	\right) \: , \;\:\;\: i = 1 \ldots n
\end{equation}
Here and in what follows $\approx$ denotes equality up to order $\varepsilon$. Therefore:
\begin{equation}
\label{disder}
	x_i(t+ \varepsilon ) \approx x_i(t) + x_i(t) 
	\left( \varepsilon \lambda ^*_{i} + \sum_{j=1}^{m} \varepsilon A^*_{ij} 
	\prod_{k=1}^{n}x_{k}^{B_{jk}}(t) \right) \: , \;\:\;\: i = 1 \ldots n
\end{equation}
If we now rescale the time so that it takes integer values, which is the usual mapping convention, we immediately find:
\begin{equation}
\label{euler}
	x_i(p+ 1) \approx x_i(p) + x_i(p) 
	\left( \lambda _{i} + \sum_{j=1}^{m} A_{ij} 
	\prod_{k=1}^{n}x_{k}^{B_{jk}}(p) \right) \: , \;\:\;\: i = 1 \ldots n
\end{equation}
In (\ref{euler}), $\lambda _{i}$ and $A_{ij}$ have been defined as $\lambda _{i} = \varepsilon \lambda ^*_{i}$ and $A_{ij} = \varepsilon A^*_{ij}$, which are dimensionless coefficients having a simple interpretation in terms of the coefficient corresponding to the interaction during a time step $\varepsilon$. 

We can now consider a slightly different possibility. The procedure is the same until 
(\ref{disder}). But then we can write, equivalently to order $\varepsilon$:
\begin{equation}
\label{disder2}
	x_i(t+ \varepsilon ) \approx x_i(t) \exp
	\left( \varepsilon \lambda ^*_{i} + \sum_{j=1}^{m} \varepsilon A^*_{ij} 
	\prod_{k=1}^{n}x_{k}^{B_{jk}}(t) \right) \: , \;\:\;\: i = 1 \ldots n
\end{equation}
We can now make the same time rescaling as in the Euler case and arrive to:
\begin{equation}
\label{qpdiscr}
	x_i(p+ 1) \approx x_i(p) \exp	\left( \lambda _{i} + \sum_{j=1}^{m} A_{ij} 
	\prod_{k=1}^{n}x_{k}^{B_{jk}}(p) \right) \: , \;\:\;\: i = 1 \ldots n
\end{equation}
As in the Euler discretization, in (\ref{qpdiscr}) we have again identified $\lambda _{i} = \varepsilon \lambda ^*_{i}$ and $A_{ij} = \varepsilon A^*_{ij}$. Consequently, the interpretation of $\lambda _{i}$ and $A_{ij}$ remains exactly the same. 

We can therefore introduce the following definition:

\mbox{}

\begin{df}
\label{df7}
{\rm For every QP differential system (\ref{glvdis}) and time step $\varepsilon$, we define its Euler discretization as the mapping 
\begin{equation}
\label{euler2}
	x_i(p+ 1) = x_i(p) + x_i(p) 
	\left( \lambda _{i} + \sum_{j=1}^{m} A_{ij} 
	\prod_{k=1}^{n}x_{k}^{B_{jk}}(p) \right) \: , \;\:\;\: i = 1 \ldots n
\end{equation}
and its QP discretization as the QP mapping
\begin{equation}
\label{qpdiscr2}
	x_i(p+ 1) = x_i(p) \exp	\left( \lambda _{i} + \sum_{j=1}^{m} A_{ij} 
	\prod_{k=1}^{n}x_{k}^{B_{jk}}(p) \right) \: , \;\:\;\: i = 1 \ldots n
\end{equation}
with $\lambda = \varepsilon \lambda ^*$ and $A = \varepsilon A^*$in both (\ref{euler2}) and (\ref{qpdiscr2}).
}
\end{df}

\mbox{}

From the previous derivation we can state a first consequence for the comparison of both discretizations:

\mbox{}

\begin{co}
\label{dco2}
{\rm Consider a QP differential system (\ref{glvdis}) and a time step $\varepsilon$. Then, 
$x_{E,i}(p) \approx x_{QP,i}(p)$ for all $i=1, \ldots , n$ and for every time $p$, where 
$x_{E,i}(p)$ and $x_{QP,i}(p)$ denote, respectively, its Euler discretization (\ref{euler2}) and its QP discretization 
(\ref{qpdiscr2}).
}
\end{co}

\mbox{}

It is then natural that certain important properties are the same for both discretizations:

\mbox{}

\begin{lm}
\label{lm3}
{\rm Consider a QP differential system (\ref{glvdis}) and its Euler and QP discretizations for a given time step $\varepsilon$, (\ref{euler2}) and (\ref{qpdiscr2}), respectively. Then
\begin{description}
\item[{\em i)\/}] A point {\bf x$_0$} $\in$ int$\{ I \!\! R^n_+\}$ is a steady state of 
(\ref{glvdis}) if and only if it is a steady state of (\ref{euler2}), and if and only if it is a steady state of (\ref{qpdiscr2}).
\item[{\em ii)\/}] If a point {\bf x$_0$} $\in$ int$\{ I \!\! R^n_+\}$ is a steady state of 
(\ref{euler2}) and (\ref{qpdiscr2}), then both discretizations (\ref{euler2}) and 
(\ref{qpdiscr2}) have the same Jacobian at {\bf x$_0$}.
\end{description}
}
\end{lm}

\begin{flushleft}
{\em Proof.}
\end{flushleft}

The proof of {\em (i)\/} is immediate. For {\em (ii)\/} it can be verified after evaluating the Jacobians while taking into account the steady state condition for each of the discretizations.$\Box$

\mbox{}

Notice that part {\em (ii)\/} of the previous lemma only states the equality of the Jacobian for both discretizations, but nothing is said about the original differential QP system. The reason is that, of course, a discretization procedure may affect the stability properties of a fixed point or, generally speaking, the local properties of the differential system. 

In fact, this may also be the case for important global properties such as, for instance, the one demonstrated in Lemma \ref{lm1}, regarding the positiveness of the trajectories of every QP mapping. In this case, QP mappings do have the expected behaviour that we intend. However, now there is an important difference between the Euler and the QP discretizations, which is that the same feature does not necessarily hold for the Euler discretization, as it can be easily verified. Therefore we thus arrive at the first example of a property which is important both from the mathematical and the applied points of view, namely the invariance of the positive orthant, which is valid for the QP discretization but not for the Euler one. This is certainly an argument in favour of the QP discretization as a more natural one for differential QP systems. 

In what follows we shall discuss another key property which is preserved by the QP discretization but not by the Euler one. 

\mbox{}

\begin{flushleft}
{\em 3.3 Commutativity of the QP discretization}
\end{flushleft}

We shall now consider the features that make the QP discretization the more consistent, from the algebraic point of view, for QP differential systems. For this, consider two 
$n$-dimensional QP systems, both having $m$ QMs and belonging to the same equivalence class. Let $\{ M^*,B \}$ and $\{ M^{\prime *},B^{\prime } \}$ be their associated matrices, respectively, and let $C$ the matrix of the QM transformation relating both systems, i.e. 
$M^{\prime *}=C^{-1} \cdot M^*$ and $B^{\prime } = B \cdot C$.

We now regard the corresponding QP discretizations of the QP differential systems for a time step $\varepsilon$. We already know that they are the QP mappings of matrices $\{ M, B\} = 
\{ \varepsilon M^*, B \}$ and $\{ M^{\prime }, B^{\prime } \} = \{ \varepsilon M^{\prime *},
B^{\prime } \}$, respectively. 

But it is then clear that the two QP mappings are related by the same QM transformation (of matrix $C$) as that which relates the original QP differential systems, since for the mapping matrices we also have $M^{\prime }=C^{-1} \cdot M$ and $B^{\prime } = B \cdot C$. In other words, we have the following:

\mbox{}

\begin{lm}
\label{lm4}
{\rm For every QP differential system, the two operations of QM transformation and QP discretization are commutative.
}
\end{lm}

\mbox{}

This property gives the QP discretization a very robust algebraic and dynamical character since, in fact, QM transformations relate topologically equivalent systems in the differential case, and topologically conjugate mappings in the discrete context. Another way to see this is by noting that the QP discretization actually operates over the entire classes of equivalence, as a bijective application associating one QP mapping to every QP differential system and vice-versa. Of course, this is the best possible situation from the point of view of the QP formalism. In addition note that: 

\mbox{}

\begin{co}
\label{dco3}
{\rm Consider a QP differential class of equivalence of matrix invariant $B \cdot M^*$, and let a QP discretization of time step $\varepsilon$. Then the resulting QP mappings form a QP class of equivalence of matrix invariant $\varepsilon B \cdot M^*$.
}
\end{co}

\mbox{}

It is not difficult to see that the same kind of considerations about commutativity apply in the case of embeddings. We omit the treatment here for the sake of brevity. However, an important consequence should be explicitly given:

\mbox{}

\begin{co}
\label{dco4}
{\rm The two operations of QP discretization and transformation to the LV canonical form are commutative.
}
\end{co}

\mbox{}

These results show that the QP discretization should be regarded as the natural one in the framework of differential QP systems. This impression is reinforced in view of the numerous common algebraic and dynamical properties and analogies shared by differential and discrete LV systems, a subject well-known in the literature \cite{hhj1,lw1}. 

In addition, it must be said at this stage that all the previous considerations could alternatively have been set in terms not of a discretization of differential QP systems, but contrarily as a continuous limit of QP mappings. In this way, it can be seen that QP differential equations can be obtained as a continuous limit of QP mappings. Both approaches are valid and completely equivalent, thus showing the symmetry between differential and discrete QP systems. In this work we have chosen to use the discretization procedure because it has some advantages in terms of easiness of exposition. The continuous limit case is thus omitted for the sake of conciseness, and because it does not lead to additional conclusions than the ones already presented. 

The results developed in this subsection are completed by considering the complementary issue of the non-existence of additional discretization procedures verifying properties to some extent analogous beyond those of the QP discretization. In this sense, we have considered mappings of the forms
\begin{equation}
\label{discrmlt}
	x_i(p+ 1) = x_i(p) \varphi_i \left( \lambda _{i} + \sum_{j=1}^{m} A_{ij} 
	\prod_{k=1}^{n}x_{k}^{B_{jk}}(p) \right) \: , \;\:\;\: i = 1 \ldots n
\end{equation}
\begin{equation}
\label{discradd}
	x_i(p+ 1) = x_i(p) + \varphi_i \left( \lambda _{i} + \sum_{j=1}^{m} A_{ij} 
	\prod_{k=1}^{n}x_{k}^{B_{jk}}(p) \right) \: , \;\:\;\: i = 1 \ldots n
\end{equation}
where the $\varphi_i$ are differentiable functions. Of course, the demonstration of non-existence results is usually a nontrivial task. Nevertheless, a careful examination of the transformation properties of mappings (\ref{discrmlt}) and (\ref{discradd}) seems to indicate that only the QP discretization has the features we are interested in. More precisely, we can state the following:
 
\mbox{}

\begin{cj}
\label{cj1}
{\rm For every QP differential system, the two operations of QM transformation and  discretization are commutative:
\begin{description}
\item[{\em i)\/}]  If and only if $\varphi_i(\xi)=a^{\xi}$, $a \in I \!\! R$, for all $i=1, \ldots ,n$, in the case of a discretization of the form (\ref{discrmlt}).
\item[{\em ii)\/}] Never, in the case of a discretization of the form (\ref{discradd}).
\end{description}
}
\end{cj}

\mbox{}

Therefore, since $a^{\xi}= e^{\xi \ln a}$, QP mappings represent the (seeming) complete set of mappings possessing the right QP formalism transformation properties. The fact that QM transformations and discretizations do not commute in the rest of the cases is sufficient in order to ensure that those discretizations do not fit into the desired algebraic scheme.

\mbox{}

\mbox{}

\begin{flushleft}
{\bf 4. Final remarks}
\end{flushleft}

The two main purposes of this work have been: 
\begin{description}
\item[{\em i)\/}] To define for the first time and describe in detail the algebraic properties of QP mappings, for which a complete discrete version of the QP formalism has been constructed. 
\item[{\em ii)\/}] To demonstrate that the algebraic properties of QP mappings are consistent from the point of view of their differential counterpart and are, in fact, similar to them. 
\end{description}

In this sense, the algebraic framework set by QP mappings and their transformations seems to be exceedingly robust due to its transformation properties and also because of the dynamical equivalence preserved by those transformations. In this sense, the definition of QP mappings opens a novel scenario not only for the establishment of algebraic or algorithmic properties, as we have already seen, but also from the point of view of applications, in which a whole range of developments regarding the possibility of robust modelling of non-linear systems can be anticipated. This impression is certainly reinforced in view of the privileged status that LV mappings play both in the QP formalism just presented and in the domain of cross-disciplinary applications, as mentioned in the Introduction. The QP discrete systems discussed so far constitue a wide, yet unexplored, generalization of those LV mappings that have led to wide developments in Physics, Population Dynamics, Chemistry and Economy. The analysis of these new possibilities shall be the subject of future research.

\mbox{}

\mbox{}

\begin{flushleft}
{\bf Acknowledgements}
\end{flushleft}

This research has been supported by a Marie Curie Fellowship of the European Community programme {\em ``Improving Human Research Potential and the Socio-economic Knowledge Base'' \/} under contract number HPMFCT-2000-00421.

\pagebreak

\end{document}